\documentclass[10pt]{article}

\usepackage{amssymb}
\usepackage{amsmath}
\usepackage{amsthm}
\usepackage{dsfont}
\usepackage{tikz}
\usepackage[round]{natbib}
\usepackage{hyperref}

\usepackage[utf8]{inputenc}
\usepackage{graphicx}
\usepackage{epsfig}
\newcommand{\R}{\mathds{R}}

\newcommand{\Z}{\mathds{Z}}

\title{Teaching Mathematics for Economists}

\author{Eric Hillebrand\footnote{Department of Economics and Business Economics, Aarhus University, Denmark, \texttt{ehillebrand@econ.au.dk}. 
}}

\date{\today}

\begin{document}

\maketitle

\abstract{In this chapter, I discuss teaching mathematical tools specifically tailored for economics students. A typical one-semester course in this area seeks to blend a range of topics: from foundational elements of subjects such as linear algebra and multivariate calculus to intermediate areas like real and convex analysis and further into advanced topics such as dynamic optimization in both continuous and discrete time. This breadth of coverage corresponds to material usually spread across multiple years in traditional mathematics programs. Given the comprehensive nature of these courses, careful selection of topics is essential, balancing numerous trade-offs. I discuss potential course sequences and instructional design choices. I then focus on conceptualizing and explaining mathematical modeling in economics. I reflect on three years of teaching an advanced undergraduate course in mathematical methods online. The latter part of the chapter offers examples and visualizations I have found particularly beneficial for imparting intuition to economics students. They cover a range of topics at different degrees of difficulty and are meant as a resource for instructors in Mathematics for Economists. Among these, I use the Ramsey model as a recurring example, especially relevant when designing a mathematical tools course with an orientation towards preparing students for macroeconomic analysis.}

\tableofcontents

\newpage

\section{The task, an ideal course sequence, and the ground truth}

From the point of view of a mathematician, teaching mathematics to economists is a nearly impossible task. At the same time, it is also fun to teach mathematics in economics programs. Students who have decided to study economics, and not mathematics, should ideally have skills in univariate calculus, linear algebra, and basic concepts of multivariate calculus (chain rule, implicit function theorem) comparable to those of students of mathematics or physics. They should have a solid grasp of ordinary differential equations, at least from a practical, solution-oriented perspective. When this is accomplished, elements of convex analysis are necessary, mainly but not exclusively for microeconomic theory, and elements of dynamic optimization and control in continuous and discrete time, particularly for macroeconomic theory. The need-to-haves are, therefore, a curious blend of elements that are at the same time very basic and quite advanced. The list of nice-to-haves is long: real analysis and probability, mathematical statistics (for econometrics), partial differential equations (for finance, mainly), stochastic optimal control (finance, macroeconomics), elements of topology, causal calculus (microeconomics), numerical optimization (all fields). I am sure that many economists will find it easy to extend this list. In addition, there are a number of areas in economics where the literature is rigorous enough that graduate students setting out to study these areas should have been introduced to rigorous proof work.

A young mathematician tasked with teaching economics students, remembering their own training where linear algebra was covered in two full semester courses, analysis in three semesters, and ordinary differential equations in one, may therefore think that one semester of linear algebra, one semester covering both univariate and multivariate calculus (ideally after the linear algebra course, so that the language of linear algebra can be used), and one semester with intermediate topics such as differential equations, elements of dynamic optimization, and convex analysis may be a good starting point at the undergraduate level. At the graduate level, one course covering analysis in a rigorous manner, for example based on Rudin's book \citep{Rudin1976}, and some advanced subjects in linear algebra, for example based on Strang's book \citep{Strang2023}, would provide a basis and a common denominator for students coming from many different places, and then, in the second year, a choice between two courses, one in dynamic analysis (differential equations and dynamic optimization) and one in convex analysis and fixed-point theory could cater to the more macro- and micro-oriented students, respectively.\footnote{Convex analysis and fixed-point theory are also employed to show existence and uniqueness of value functions in dynamic optimization, and there are mathematics-for-economists textbooks that cover this, for example, \cite{Corbae2009}.}

In the morning, our mathematician runs into the econ department head on the way from the parking lot and mentions his rough outline for a course sequence. With a puzzled look, the department head says, ``That sounds intriguing. We thought you would be the perfect fit to teach the Chiang course.'' Before our mathematician has any chance to rephrase the thought ``What on earth is a Chiang course?'' the department head continues, ``Well, we do have two math courses on the books, because we're quantitatively quite strong, one at the undergrad and one at the grad level! If you'd like to revamp a course or design a new one, of course you're welcome to write a proposal.'' 

Later that day, the senior econometrician, a very kind and smart man with a knack for mentoring, whom our young hero immediately understood it was important to listen to, drops by. ``Hey, I heard you pitched some math courses to the department head. It's a great thought to do more math, and I would love to see it, but it's going to be an uphill battle to get more courses in than we have right now. And leave the matrix algebra for the graduate level; our undergraduates cannot stomach it. Believe me. I tried. Same goes for complex numbers, and we don't really need them in economics, anyway.'' \nocite{Chiang1984}

Thus goes up in smoke our mathematician’s dream of sharing with their students their own undergraduate epiphany that once you allow for higher dimensions, for example, 2$\times$2, the need for $i$ disappears, since the equation
\[
X^2 = -I
\]
has the perfectly valid integer solution
\[
X = \left[\begin{array}{cc}
0 & -1 \\ 1 & 0
\end{array}\right]\in\Z^{2\times 2}.
\]
Lo and behold, the eigenvalues of this matrix are $i$ and $-i$, but there's really no need to insist on computing them. 

The view that linear algebra cannot be communicated to undergraduate economics students is very common. It pervades the design of many textbooks, even otherwise excellent ones, such as \cite{Sydsaeter2022} and \cite{Sydsaeter2008}.\footnote{A notable exception is \cite{Simon1994}.} I think that this view is one of the major hindrances to a better math education within economics. It drastically reduces the speed with which progress can be made towards the material that is necessary for understanding and contributing to the economics literature---or that would enable the economics student to appreciate when neighboring disciplines, such as engineering, talk about exactly the same things as economists do \citep{Luenberger1969}. The situation results in many serious students, who are ambitious about an academic career in economics, resorting to taking courses at mathematics departments. On the other hand, maybe it is also better that way. There is a certain mathematical culture that is very useful for an economist working with rigorous material to have experienced, and it is probably best when it comes from the horse's mouth. The serious student quickly realizes that even the ostensibly basic material in linear algebra, analysis, and differential equations presents itself at a completely different depth and brilliance when studied carefully. If they made the commendable decision to take some of the undergraduate courses at the mathematics department, for example, as a first-year economics Ph.D. student, they are then often met with the---in this case---absurd verdict that they cannot get graduate credit for taking undergraduate classes. This is, in fact, often the motivation for having a graduate-level mathematics-for-economics class on the books for the economics department in the first place. The treatment of multivariate and convex analysis in these courses often does not go beyond undergraduate material for mathematics students. Dynamic optimization is different, though.

Returning to linear algebra, should it really be beyond the capacity for abstraction of a 19- to 23-year-old student of economics that a vector is a point in a properly defined space that behaves under addition with other vectors and multiplication with scalars as you would expect? That multiplication of a vector by a matrix, for example the product $y=f(x)=Ax$, where
\[
A = \frac{1}{2}\left[\begin{array}{cc}
5 & -1 \\
-1 & 5
\end{array}\right],\qquad 
x = \left[\begin{array}{c}
1 \\ 3
\end{array}\right],
\]
is given by the two lines of addition and multiplication operations
\[
\left[\begin{array}{c}
\frac{5}{2} \cdot 1 - \frac{1}{2}\cdot 3 \\
-\frac{1}{2} \cdot 1 + \frac{5}{2}\cdot 3
\end{array}\right]?
\]
That every multiplication of a vector $x$ by a matrix $A$ amounts to a rotation and a scaling of $x$?
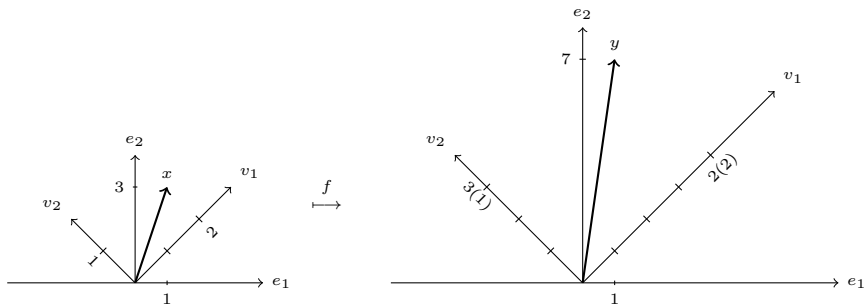
\begin{center}
\begin{tikzpicture}[scale=0.85]
\draw[->] (-2,0) -- (2,0) node[anchor=west] {$\scriptstyle e_1$};

\draw (.5,1pt) -- (.5,-1pt) node[anchor=north] {$\scriptstyle 1$};

\draw[->] (0,0) -- (0,2) node[anchor=south] {$\scriptstyle e_2$};

\draw (1pt,1.5) -- (-1pt,1.5) node[anchor=east] {$\scriptstyle 3$};

\draw[->,thick] (0,0) -- (.5,1.5) node[anchor=south] {$\scriptstyle x$};

\draw[->] (0,0) -- (1.5,1.5) node[anchor=south west] {$\scriptstyle v_1$};

\draw[->] (0,0) -- (-1,1) node[anchor=south east] {$\scriptstyle v_2$};

\draw (.45,.55) -- (.55,0.45);

\draw (.95,1.05) -- (1.05,.95);

\draw (1.2,.8) node[rotate=45] {$\scriptstyle 2$};

\draw (-.55,.45) -- (-.45,.55);

\draw (-.65,.35) node[rotate=-45] {$\scriptstyle 1$};

\draw (3,1.5) node {$\scriptstyle f$};

\draw (3,1.2) node {$\scriptstyle \longmapsto$};

\draw[->] (4,0) -- (11,0) node[anchor=west] {$\scriptstyle e_1$};

\draw (7.5,1pt) -- (7.5,-1pt) node[anchor=north] {$\scriptstyle 1$};

\draw[->] (7,0) -- (7,4) node[anchor=south] {$\scriptstyle e_2$};

\draw (7.05,3.5) -- (6.95,3.5) node[anchor=east] {$\scriptstyle 7$};

\draw[->,thick] (7,0) -- (7.5,3.5) node[anchor=south] {$\scriptstyle y$};

\draw[->] (7,0) -- (10,3) node[anchor=south west] {$\scriptstyle v_1$};

\draw[->] (7,0) -- (5,2) node[anchor=south east] {$\scriptstyle v_2$};

\draw (7.45,.55) -- (7.55,0.45);

\draw (7.95,1.05) -- (8.05,.95);

\draw (8.45,1.55) -- (8.55,1.45);

\draw (8.95,2.05) -- (9.05,1.95);

\draw (9.2,1.8) node[rotate=45] {$\scriptstyle 2(2)$};

\draw (6.45,.45) -- (6.55,.55);

\draw (5.95,.95) -- (6.05,1.05);

\draw (5.45,1.45) -- (5.55,1.55);

\draw (5.35,1.35) node[rotate=-45] {$\scriptstyle 3(1)$};

\end{tikzpicture}
\end{center}
That the angle of rotation and the scaling factor are, however, in general not independent of $x$, that they can vary from vector to vector?\footnote{That only specific matrices always have the same rotational angle and the same scaling independent of $x$? These matrices are $A=\left[\begin{array}{cc}
\alpha \cos\theta & -\alpha\sin\theta \\
\alpha \sin\theta & \alpha\cos\theta
\end{array}\right]$, $\alpha\in\R$ scaling factor, $\theta\in [0, 2\pi]$ rotational angle.}
That despite this dependence of the scaling and rotation on the argument vector of the function, many matrices allow for a transformation of the coordinate system such that the scaling in that transformed coordinate system \textit{along those axes} is always the same? In our example, the transformed coordinate system is given by the vectors
\[
\left(\left[\begin{array}{c} 1 \\ 1\end{array}\right],\, \left[\begin{array}{c} -1 \\ 1\end{array}\right]\right).
\]
In the multiplication $y=Ax$, first the old coordinates $(1,\, 3)$ of
$x$ are switched to the new coordinates $(2,\, 1)$ with respect to this system.  Then, the new coordinates are stretched by the eigenvalues, the first coordinate by the first eigenvalue, the second coordinate by the second eigenvalue, to $(4,\, 3)$.  These are the new coordinates of $y$. Finally, the
new coordinates are converted back to old coordinates $(1,\, 7)$.
\begin{align*}
y = Ax &= \underbrace{\left[\begin{array}{cc} 1 & -1 \\ 1 & 1\end{array}\right]}_{P} \underbrace{\left[\begin{array}{cc} 2 & 0 \\ 0 & 3\end{array}\right]}_{\Lambda} \underbrace{\frac{1}{2} \left[\begin{array}{cc} 1 & 1 \\ -1 & 1\end{array}\right]}_{P^{-1}} \left[\begin{array}{c}
1 \\ 3
\end{array}\right] \\
&= \underbrace{\left[\begin{array}{cc} 1 & -1 \\ 1 & 1\end{array}\right]}_{P} \underbrace{\left[\begin{array}{cc} 2 & 0 \\ 0 & 3\end{array}\right]}_{\Lambda} \left[\begin{array}{c}
2 \\ 1
\end{array}\right]\\
&= \underbrace{\left[\begin{array}{cc} 1 & -1 \\ 1 & 1\end{array}\right]}_{P}  \left[\begin{array}{c}
4 \\ 3
\end{array}\right] = \left[\begin{array}{c}
1 \\ 7
\end{array}\right]\\
\end{align*}
I'm convinced that this is within the geometric and algebraic reach of an undergraduate economics student. They may complain about the method for finding the eigenvalue decomposition, but it involves finding the roots of a polynomial and solving a system of equations, all things they have seen before in some way.

Be that as it may, our young mathematician now has to make do with the one-semester undergraduate and the one-semester graduate course plus the senior econometrician's advice. At the undergrad level, they will end up teaching some version of multivariate calculus with two or at most three variables, writing scalar partial derivatives, not gradients or Jacobians. They will likely do some version of optimization under constraints but cannot really talk about constraint qualification. Some review of integration of functions of one variable will likely be done.  Perhaps some scalar first-order ordinary differential equations are covered, certainly not systems or higher-order equations, where many deeply satisfying insights into the relations of ordinary differential equations with linear algebra await. If the young mathematician is very organized and does not follow their impulse to delve too deeply into the subjects, some dynamic optimization can be covered, with some luck, before the semester runs out. This would be either in continuous or discrete time---most likely not both---and, of course, scalar only, probably deterministic, not stochastic. The students will still find the course dense and challenging. The chalkboard is full of math, after all. 

Most mathematics-for-economists textbooks represent some compromise along these lines. Some opt for the voluminous or even two-volume approach \citep{Simon1994,Hoy2001,Klein2002,Sydsaeter2022}, others opt for the small-and-lean approach \citep{Weintraub1982,Novshek1993}. In the graduate course, remedial linear algebra now needs to be taught. This easily takes a third of the semester. There's still time to teach some real analysis, proof techniques, and advanced topics in convex analysis and dynamic optimization. Similar to the undergraduate situation, a balance has to be struck, and the specific solution usually ends up depending on a combination of the popular specializations that Ph.D. students take in the program at hand and the personal preferences of the instructor. Comparing graduate texts, it quickly becomes clear that the subject can be treated ``right'' in many different ways, see, for example \cite{Ok2007}, tailored towards econometrics, \cite{Corbae2009}, tailored towards macroeconomics, and \cite{Sydsaeter2008}, trying to balance macro- and microeconomics, leaning towards dynamic optimization. Both at undergrad- and grad-level, single one-semester courses prevail, which have some umbrella title such as ``Mathematics for Economics'', ``Mathematical Methods in Economic Analysis'' or similar, and do not branch out into subjects or purposes.  Together with the common perception that textbooks should be relatively new, even though the subject material taught has essentially not changed in over 50 years, this makes it difficult to assign a number of wonderful monographs with more specific titles \citep[for example,][and this is certainly not an exhaustive list]{Seierstad1987,Dixit1990,Intriligator2002,Dhrymes2013}.\footnote{M.D. Intriligator noted in the 2002 SIAM Classics edition of his 1971 ``Mathematical Optimization and Economic Theory'' that the basic methods covered remain the same as in 1971, and that the only change he would have made was putting more emphasis on duality, reflecting developments in mathematical economics.}

\section{Mathematical modeling in economics}

This essay starts from the assumption that students of economics should be taught a level of mathematics that enables them to comprehend (undergraduate and master's level) and contribute to (Ph.D.) parts of the academic economics literature. The ideal course sequence proposed above operates from this basic assumption.\footnote{The subject and the courses discussed here are ``mathematics for economists/economics'', i.e., tools courses, and not ``mathematical economics.'' There is no universally agreed upon distinction, and there is a good deal of confusion, but it makes sense to me to reserve the latter term, as many economists do, for formal economic equilibrium theory, represented by texts such as \cite{Debreu1959}, \cite{Arrow1963}, \cite{Arrow1981,Arrow1982,Arrow1986}, or, for example, \cite{Stigum1990}. Mathematical economics progresses, and a course on this subject likely will need more frequent updating than a tools course.} Whether the mathematical sophistication employed in the economics literature is always useful for good economics is a different question that reasonable people can argue about, and I do not address that question here. Inevitably, at some point, the math-for-econ instructor is confronted with the student's question of what all these mathematical models are good for, if the instructor ``believes'' in them, if they are ``true'', or some variation on that theme. This is a good question that deserves encouragement.

Probably the most famous answer to this question is George E.P. Box's \citep{BOX1979201}:

\begin{quote}\it
``All models are wrong but some are useful,'' 
\end{quote}
with an entire Wikipedia article devoted to it.\footnote{\url{https://en.wikipedia.org/wiki/All_models_are_wrong}}

The quote pointedly describes the trade-off between accuracy and usefulness of a model. In my own lectures, I like to illustrate this trade-off with two extreme examples: At one extreme, imagine mythological explanations for why one should seek higher ground when the water at the beach recedes.\footnote{\url{http://news.bbc.co.uk/2/hi/south_asia/4181855.stm}} The mental model used here is almost certainly entirely inaccurate, but tremendously useful. I'd certainly prefer surviving with it over drowning while relying on  differential equations, because I, again, forgot the integration constant. At the other extreme, consider a perfectly accurate model, such as a map of a city at scale 1:1 and in three dimensions \citep{Borges1999,Casati2024}. One is as lost in the map as in the actual city, so it is entirely useless.\footnote{I think that we are currently observing a perfectly accurate model of money in cryptocurrencies. These can serve all three functions of money, at least in theory, and therefore satisfy the standard textbook definition. However, another and unrelated definition of money that one often hears when talking with bankers is that money is a liability of a money-issuing, usually central bank, i.e., that there is a creditworthy institution that has the money issue on the passive side of the balance sheet. If the failure of cryptocurrencies to satisfy this definition becomes problematic, this will at the same time demonstrate the uselessness of cryptocurrencies as a model for money, because they do not improve our understanding of money.} There are other dimensions that play into the trade-off, such as the state of understanding of the causal processes at work, the effort required to run the model, and the purpose of the model. 

John von Neumann defined science solely as model building \citep{Neumann1995}:
\begin{quote}\it
``The sciences do not try to explain, they hardly even try to interpret, they mainly make models. By a model is meant a mathematical construct which, with the addition of certain verbal interpretations, describes observed phenomena. The justification of such a mathematical construct is solely and precisely that it is expected to work -- that is, correctly to describe phenomena from a reasonably wide area.''
\end{quote}
This seems overly restrictive to me, but it focuses on interesting aspects of model building: Explanations and interpretations of causal chains are usually not the result of models; they are premises of the model building cycle. The model is a metaphor of the causal processes at work according to current understanding, and it helps organize thinking and put features of interest into sharp relief. This aids communication of ideas between researchers (and between teachers and students), and it helps to pose sharp questions. If the model can predict, great, but this is often not the main interest in the model. It is more important that the model can describe the phenomenon of interest with minimal assumptions. In economics, this minimum is often referred to as ``from first principles.'' Examples of first principles are diminishing marginal returns, concave or quasi-concave utility functions, time preference, and so on. Modelers find it satisfying if an empirical regularity can be recreated in the model with a minimal set of assumptions. Contrary to von Neumann's dictum, they then like to say that the model ``explains'' the regularity. 

In statistical and econometric models, the empirical reality to be captured is the data. In structural models, it is again a minimal set of assumptions borrowed from economic theory that are seen as desirable. The role of the random variables in the econometric model is conceptualized on a spectrum: From \textit{``\ldots errors appear because the model is a simplification of reality [and may omit variables, EH] \ldots error is associated with the collection and measurement of the data \ldots''} \citep[p. 58]{Pindyck1998}, a perspective in which noise is added, to \textit{``\ldots the structural shocks [\ldots] having a distinct economic interpretation. [\ldots] Current observations of the data may be viewed as a weighted average of current and past structural shocks.''} \citep[p. 2,3]{Kilian2017}, a perspective in which the errors are the primitive, and the structural relations are only the description of the channels through which the errors propagate. 

The litmus test is whether the model can produce output that matches the distributional characteristics of the data. This can be assessed by checking whether the residuals of the estimated model satisfy the distributional assumptions on the error terms according to diagnostic statistics. Alternatively, it can be checked by examining whether the moments of the model output match the moments of the data. In reduced-form models, it is minimal complexity, or parsimony, that is often sought. Parsimony always refers to a minimum number of parameters, in some sense. They can be associated with explanatory variables, lags of variables, deterministic or stochastic features of the model, etc. Again, model output consistent with the statistical properties of the data is aimed for. Never is the model conceptualized as an accurate description of reality, not even when the ``nonlinearities'' are taken into account. To take the model literally is to not take it seriously \citep{Thompson2022}. The choice of which features of reality to emphasize in a model and which ones to tone down depends on a mixture of all the aspects mentioned above, plus, to no small extent, the skill and experience of the modeler. It is often said that modeling is an art rather than a science. Modeling is certainly a creative process.

As a caveat, it is my understanding that philosophy of science has come a long way since the people whom I'm citing here were active and helped shape the epistemological background of the training of my generation \citep{Schindler2018}. On a lighter note, modeling is taken seriously in economics to the point of inviting self-parody. I recommend the unparalleled ``Life Among the Econ'' \citep{Leijonhufvud1973}:
\begin{quote}\it
``The Econ word for caste is `field.' Caste is extremely important to the self-image and sense of identity of the Econ, and the adult male meeting a stranger will always introduce himself with the phrase `Such-and-such is my field.' [\ldots] A comparison of status relationships in the different `fields' shows a definite common pattern. The dominant feature, which makes status relations among the Econ of unique interest to the serious student, is the way that status is tied to the manufacture of certain types of implements, called `modls.' The status of the adult male is determined by his skill at making the `modl' of his `field.' The facts (a) that the Econ are highly status-motivated, (b) that status is only to be achieved by making `modls,' and (c) that most of these `modls' seem to be of little or no practical use probably accounts for the backwardness and abject cultural poverty of the tribe.''
\end{quote}

\section{Teaching Mathematics for Economists online}

At Aarhus University, at the time of writing, we are in the relatively luxurious position of having two-and-a-half semester courses in mathematics for economists, one-and-a-half at the undergraduate level and one at the graduate level. The half-semester course is the first half of ``Principles of Mathematics and Statistics,'' taught in Danish in the first semester, and the mathematics part is based on selected chapters from the book by \cite{Sydsaeter2022}. There are about 200 students per year. The full semester course is called ``Mathematics for Economists,'' which I had renamed in a long process from ``Mathematical Economics'' after I took over the course in 2014. It is taught in Danish in the third undergraduate semester, and it is based on selected chapters from the book by \cite{Sydsaeter2008}. There are about 160 students. In this course, I cover linear algebra in three weeks, multivariate calculus in the next three weeks, static optimization (global, local, and constrained extrema) in one-and-a-half weeks, integration of functions of several variables (one week), first-order ordinary differential equations (one week), second-order linear differential equations and two-dimensional systems (one week), control theory (two weeks), and difference equations and dynamic optimization in discrete time (one-and-a-half weeks), for a total of fourteen weeks. The exam format is a 20-minute oral exam. We screen the incoming undergraduate students for mathematics; the Danish high-school system allows students to take math at different levels, and we require the highest.

The Ph.D. course was taught for years in English by a colleague from the mathematics department, who used \cite{Rudin1976}. There are about ten students per year. After his retirement, the course has been taught by different colleagues who have worked with various texts \citep{Boyd2004,Johnsonbaugh2010,Ok2007,Resnick2014}. The emphasis on real analysis represents the least common denominator for all fields, since it directs the focus towards rigorous proof training. The probability text was used supplementarily. We have a number of international Ph.D. students, but the majority are Danish with their undergraduate degree from Aarhus, so they have already taken the third-semester undergraduate course.

This was also the situation in 2020, when COVID-19 hit. I had been toying with the idea of taking the third-semester course with about 160 students online for a long time before that. As noted above, the subject matter has not changed substantially in decades, and it certainly does not change from year to year. This makes the course a near-perfect candidate for recording videos and experimenting with an online format. I had some experience with recording videos on mathematics for economists from my YouTube channel, which I had started sometime in 2016.\footnote{\url{https://www.youtube.com/channel/UCJjiXpAzy_L6uneX5uoTihQ}} I took the situation in 2020 as an opportunity to give it a try. The easiest part was to choose a video format. I went for a three-part layout, with half of the screen filled with a projection of my tablet, which I used to write on, a quarter of the screen filled with the slides and another quarter with a headshot video of me. This way, I figured, the students would have the main lecture experience in a nutshell: the instructor, the blackboard, the slides.\footnote{An example, in Danish, can be seen here: \url{https://www.youtube.com/watch?v=HDXOZnMLaes}.}

In 2020, the lockdown situation also necessitated taking the tutorials, run by student teaching assistants, online. This was the more difficult problem to solve. In the physical format, I had been giving weekly homework sets with exercises, some from the book and some of my own, and in the in-person tutorials, the students had been discussing these exercises. For the online format, I decided to drastically reduce the number of mathematical exercises to work through and instead added an exercise in which each working group of three or four students had to prepare an eight-minute video with a presentation of the main elements of the week's material. The students recorded these videos with their smartphones, for example by placing them on a stack of books with the camera aimed at a sheet of paper. The eight minutes were chosen because this was about the time a student presented in the 20-minute exam. The rest of the time was spent with questions and answers. 

In the online tutorial, half of the time was spent watching and discussing the videos and the other half discussing the mathematical exercises. My main intention was to prepare the students for the likely scenario of online oral exams, which did in fact occur, both in 2020 and 2021, due to government-ordered lockdowns. (The exams are in January.) In preparing the videos, the students could establish a routine with their improvised recording setups and with the situation of explaining math to a camera and a microphone.\footnote{The level of students' creativity and ingenuity was something to behold. A stellar moment was a student who sat in her kitchen for the online oral exam with the smartphone camera aimed at her white refrigerator door, on which she answered our questions impeccably with a whiteboard marker. It was beautiful.} Another intention was to give the students the opportunity to see other students in exam-like presentations and to make use of the common asymmetry that one is better at judging than at performing oneself. The hope was that they would catch some of the common mistakes I had seen in the exams early on, when other students made them in the tutorials. This worked to a degree, but meaningful feedback on the videos between the different working groups remained an exception. The students clearly did not feel comfortable critically discussing their peers' videos and when asked, resorted to short and generally positive and encouraging statements (``great video!''). I got feedback myself from our Center for Educational Development that in their experience, getting students to give each other meaningful feedback requires holding workshops with them.  

I ended up running the online format for three consecutive years, from 2020 to 2022. In 2021 and 2022, the tutorials were in person. As noted above, the oral exams had to be held online in 2020 and 2021. I met with the teaching assistants once a week to discuss what was going on in the tutorials and the exercises. In the in-person format, there had been years where, to my embarrassment, I had difficulties remembering the TAs' names by the end of the semester. In the online format, this became unthinkable. You could have woken me up in the middle of the night and I could have given you the full names of my (usually four) TAs. This was due in no small part to the realization that by choosing the online format, I had transferred the main teacher role in the teacher-student relationship to them. They were the ones who met with the students once a week and who were the course's ``face'' to the students. They saw the students working, answered their questions, and only related to me the questions and issues they could not solve themselves. Their online office hours were used heavily, whereas my own online office hours were largely perfunctory. 

To maintain traction with the students and stay in touch, I visited the tutorials on a sampling basis, both in the online and in-person tutorial versions. Since there are usually eight tutorial groups of about 20 students, I could visit each at most twice per semester. I got very enthusiastic about this process. In the in-person format, I lecture in front of 160 students in a large hall, and interaction is limited to the two handfuls of students who come up in the break or after the lecture. I only see the students working, and learn their names, in the oral exam at the end of the semester. In the online format, I saw them working on problems on the blackboard from the first week on, identified the good students throughout the semester (TA material for later), and generally had the feeling of being much closer to the students. It took me a while to understand that this feeling was thoroughly asymmetric. In the in-person format, the students see me twice a week for 90 minutes, hear me talk, watch me work, and see me get chalk-dusty. Now they saw me in a small video window on YouTube and otherwise at most twice in tutorials, mostly sitting and listening and then giving feedback on presentations, likely not their own, and making some general comments, and that was it. They felt I was far away.

In the oral exams, the online students fared substantially worse. They were perhaps marginally better at presenting what they had prepared, but they were generally much worse at answering anything that lay outside their preparation. This was often true for obvious and predictable questions, such as when the topic was titled ``Second-order ordinary differential equations and systems in the plane,'' and they had prepared a presentation on second-order equations and were then surprised when I asked them about systems. This also happened in spite of explicit and coordinated warnings to this effect in the tutorials. I was puzzled about the possible causes, and I still haven't come to a satisfying conclusion. My working hypothesis is that the regular presentations in the tutorials, which were a much more relaxed situation compared to the exam and were followed only by friendly feedback from the student TA, had framed the exam as just another one of those. But I'm not sure.

Another effect I observed has also been documented in a large-scale randomized controlled trial in the literature \citep{cacault2021distance}. The good students thrived in the online format. They were intrinsically motivated and interested in engaging with the material, and so they could enjoy the flexibility the online format gave them. They liked it and were as good as ever. The not-so-good students, however, had a worse outcome than in the in-person format. They struggled to keep up with the weekly flow of material, and in informal discussions with the student TAs, a higher-than-normal percentage of students gave the feedback that they were lost. They did not show up for my online office hours, probably for fear of embarrassment, and the student TA office hours were a reaction to this pattern. Still, only a few dared come to their TAs, and even though the student TA office hours were a success, they were also much frequented by the good students who came to discuss finer points, for which, of course, there also should be time. Again, I can only hypothesize as to the possible causes. One is that upon hitting a rough patch in working with material, it is, of course,  much easier to close down a laptop than it is to leave a room full of people. But also in the in-person format most of the deeper understanding is happening at home, so, this is only a partial explanation. Another element may be that in the in-person format, a struggling student can physically see peers successfully engaging with the material, and this may give them a much-needed kick---maybe out of a sense of competitiveness, maybe because they think, “If they can do it, I should be able to as well,” maybe simply because it sets a certain tone and positive attitude. Wherever the effect originated, it was clearly discernible.

In my own time economy, after the first year, the online course took almost, but not quite, as much time as the in-person version. In the first year, recording the videos took an inordinate amount of time, and for about two months, I was not doing anything else. In the following years, visiting the tutorials, holding online office hours, engaging with the TAs at a deeper level, and fielding student emails, all amounted to just slightly less work than the in-person version. The total time commitment to students and TAs was about the same, but I did not have to prepare classes.

In the fall of 2023, I switched back to the in-person format. The reason was a combination of all of the above, but the main motivation was that I became a professor to teach and to research. I enjoy teaching. In the online format, my role was reduced to that of a teaching organizer and facilitator. I did many teaching-related tasks, but as noted above, the actual teachers were the TAs. Not having to prepare classes, I realized to my chagrin that I had become a bit rusty in the subject material. Returning to the in-person format, I eliminated the presentation exercise from the homework sets and increased the number of problems to solve---pencil to paper, or stylus to tablet. The students became better again at answering exam questions outside their preparations, or at least, they were not surprised by them. There is a time trend, however, that is pointing in the wrong direction. Occasionally, I get the comment that written exercises are not the best preparation for oral exams, because that sounds true. This is difficult to reply to without telling the long story, but the answer is, as I now know, that they actually are.  Some students ask if I could make the old videos available, which I do not think is a good idea. The video linked above is accessible solely for the purpose of this paper.

\section{Examples and exercises}

In this section, I present selected examples and exercises. They are meant as a resource for instructors in Mathematics for Economists and may be freely used. My github repository houses the Latex code for this section, including the tikz bits for the illustrations.

\subsection{The geometry of the determinant}

\begin{quote}\it
``The determinant of a square matrix $A$, denoted by $|A|$, is a uniquely defined scalar (number) associated with that matrix.'' \citep[p. 93]{Chiang1984}
\end{quote}

\smallskip

The determinant of an $n\times n$ square matrix is the $n$-dimensional volume of the parallelepiped spanned by the column (or row) vectors of that matrix. There is a sign convention for this volume, similar to sign conventions in integration. Thus, in the case of $n=2$, the determinant is the area of a parallelogram.

\subsubsection{The determinant of a 2$\times$2 matrix}

In this example, I show how the determinant of a 2$\times$2 matrix describes the area of the parallelogram spanned by the column vectors, using arguments from geometry and trigonometric functions that the students should have seen in high school. I also make use of the inner (or scalar) product, which I introduce right before determinants.

Consider the matrix
\[ A= \left[ \begin{array}{cc}
3 & 1\\
1 & 4
\end{array}\right] =: [v_1\; v_2].
\]

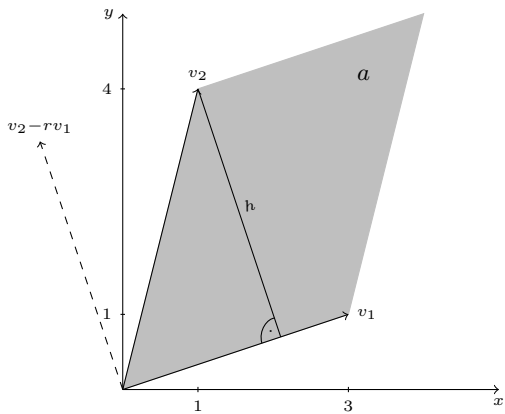
\begin{center}
\usetikzlibrary{patterns}
\begin{tikzpicture}
\filldraw[color=lightgray] (6,0) -- (7,4) -- (10,5) -- (9,1) -- (6,0);

\draw[->] (6,0) -- (11,0) node[anchor=north] {$\scriptstyle x$};

\draw[->] (6,0) -- (6,5) node[anchor=east] {$\scriptstyle y$};

\draw[->] (6,0) -- (7,4) node[anchor=south] {$\scriptstyle v_2$};

\draw[->, dashed] (6,0) -- (4.9,3.3) node[anchor=south] {$\scriptstyle v_2 - r v_1$};

\draw[->] (6,0) -- (9,1) node[anchor=west] {$\scriptstyle v_1$};

\draw (9,1pt) -- (9,-1pt) node[anchor=north] {$\scriptstyle 3$};

\draw (7,1pt) -- (7,-1pt) node[anchor=north] {$\scriptstyle 1$};

\draw (6.03,4) -- (5.97,4) node[anchor=east] {$\scriptstyle 4$};

\draw (6.03,1) -- (5.97,1) node[anchor=east] {$\scriptstyle 1$};

\draw (8.1,.7) -- (7,4);

\draw (7.85,0.617) .. controls (7.8,.79) and (7.93,.95) .. (8.02,.95);

\draw (7.97,.78) node {$\scriptstyle \cdot$};

\draw (7.5,2.4) node[anchor=west,rotate=15] {$\scriptstyle h$};

\draw (9,4) node[anchor=south west] {$a$};
\end{tikzpicture}
\end{center}

The columns $v_1$ and $v_2$ of $A$ span a parallelogram $\{r_1 v_1 + r_2
v_2 | 0\le r_{1}, r_2\le 1\}$.  To determine the area $a$ of the
parallelogram, we need to find the height $\| h\|$. We have that $\| h\| = \| v_2-r v_1\|$
for some $r\in\R$, and $v_2-r v_1$ is perpendicular to $v_1$. This results in the following condition for $r$:
\begin{align*}
\langle h, v_1\rangle=\langle v_2- r v_1, v_1\rangle &= \langle v_2, v_1\rangle - r\, \|v_1\|^2 \stackrel{!}{=} 0\\
 &= 7 - 10r \; \Longrightarrow r=\frac{7}{10},
\end{align*}
where $\langle x, y\rangle = x^T y = x_1y_1 + x_2y_2$ denotes the inner product, or scalar product, of the vectors $x$ and $y$. Then,
\begin{align*}
h=v_2 - r v_1 &= \left[\begin{array}{c} 1\\ 4 \end{array}\right] - \frac{7}{10} \left[\begin{array}{c} 3\\ 1 \end{array}\right] = \left[\begin{array}{c} -1.1\\ 3.3 \end{array}\right],\\
\| h\| &= \| v_2 - r v_1\| = \sqrt{1.1^2 + 3.3^2},\\
a &= \| v_1\|\, \| h\| = \sqrt{10} \, \sqrt{1.1^2 + 3.3^2} = 11,
\end{align*}
and we have found the area of the parallelogram to be $11$ units.

Alternatively, define
\[
v_1' := \frac{v_1}{\| v_1\|} = \frac{1}{\sqrt{10}} \left[\begin{array}{c} 3\\ 1 \end{array}\right] ,\; v_2' := \frac{v_2}{\| v_2\|} = \frac{1}{\sqrt{17}} \left[\begin{array}{c} 1\\ 4 \end{array}\right].
\]
The vectors $v_1'$ and $v_2'$ have unit length. Consider the parallelogram
spanned by $v_1'$ and $v_2'$ and denote its area $a'$.

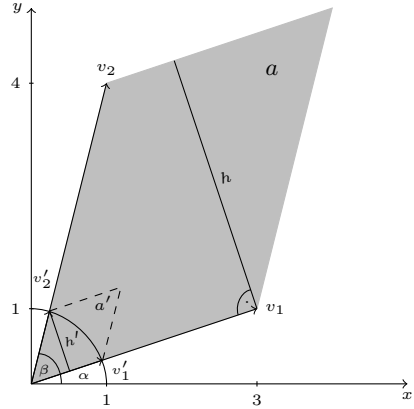
\begin{center}
\begin{tikzpicture}

\filldraw[color=lightgray] (0,0) -- (1,4) -- (4,5) -- (3,1) -- (0,0);

\draw[->] (0,0) -- (5,0) node[anchor=north] {$\scriptstyle x$};

\draw[->] (0,0) -- (0,5) node[anchor=east] {$\scriptstyle y$};

\draw[->] (0,0) -- (1,4) node[anchor=south] {$\scriptstyle v_2$};

\draw[->] (0,0) -- (3,1) node[anchor=west] {$\scriptstyle v_1$};

\draw (3,1pt) -- (3,-1pt) node[anchor=north] {$\scriptstyle 3$};

\draw (1,1pt) -- (1,-1pt) node[anchor=north] {$\scriptstyle 1$};

\draw (0.03,4) -- (-0.03,4) node[anchor=east] {$\scriptstyle 4$};

\draw (0.03,1) -- (-0.03,1) node[anchor=east] {$\scriptstyle 1$};

\draw (3,1) -- (1.9,4.3);

\draw (2.75,0.917) .. controls (2.7,1.09) and (2.83,1.25) .. (2.92,1.25);

\draw (2.87,1.08) node {$\scriptstyle \cdot$};

\draw (2.4,2.7) node[anchor=west,rotate=15] {$\scriptstyle h$};

\draw (3,4) node[anchor=south west] {$a$};

\draw (0,1) .. controls (.5,1) and (1,.5) .. (1,0);

\draw[->] (0,0) -- (.949,.316);

\draw (1.2,.18) node {$\scriptscriptstyle v_1'$};

\draw[->] (0,0) -- (.243,.970);

\draw (.15,1.4) node {$\scriptscriptstyle v_2'$};

\draw[dashed] (.243,.970) -- (1.192,1.286);

\draw[dashed] (.949,.316) -- (1.192,1.286);

\draw (.17,.17) node {$\scriptscriptstyle \beta$};

\draw (.1,.4) .. controls (.25,.4) and (.4,.25) .. (.4,0);

\draw (.7,.1) node {$\scriptscriptstyle \alpha$};

\draw (.51,.17) -- (.243,.97);

\draw (.55,.6) node[rotate=15] {$\scriptscriptstyle h'$};

\draw (1.192,1.286) node[anchor=north east] {$\scriptstyle a'$};

\end{tikzpicture}
\end{center}

The coordinates of the unit vectors
are $v_1'=(\cos\alpha,\,\sin\alpha)^T$ and $v_2'=(\cos\beta,\,
\sin\beta)^T$.  Since $v_1$ is a multiple of $v_1'$, and $v_2$ is a
multiple of $v_2'$, we can write for $\rho=\sqrt{17}$, $\sigma=\sqrt{10}$,
\begin{align*}
v_1 &= \rho v_1' = \rho \left[\begin{array}{c}\cos\alpha\\\sin\alpha\end{array}\right] = \| v_1\| v_1' = \sqrt{17}\, v_1',\\
v_2 &= \sigma v_2' = \sigma \left[\begin{array}{c}\cos\beta\\\sin\beta\end{array}\right] = \| v_2\| v_2' = \sqrt{10}\, v_2'.\\
\end{align*}
The length of $h'$ is given by
\[
\| h'\| = \frac{\| h'\|}{\| v_2'\|} = \sin(\beta-\alpha).
\]
The area $a'$ is therefore
\begin{align*}
a' &= \| h'\|\, \|v_1'\| = \| h'\| = \sin(\beta-\alpha),\\
 &= \cos\alpha \sin\beta - \cos\beta\sin\alpha,
\end{align*}

Define the matrix
\[
A' = \left[\begin{array}{cc}
a_{11}' & a_{12}' \\
a_{21}' & a_{22}'
\end{array}\right] =
\left[\begin{array}{cc}
\cos\alpha & \cos\beta \\
\sin\alpha & \sin\beta
\end{array}\right] = [v_1'\; v_2'],
\]
then we notice that
\[
a' = \cos\alpha \sin\beta - \cos\beta\sin\alpha = a_{11}' a_{22}' - a_{21}' a_{12}'.
\]
Consider the matrix $A$ that spans the larger parallelogram
\[
A = \left[\begin{array}{cc}
a_{11} & a_{12}\\
a_{21} & a_{22}
\end{array}\right] =
\left[\begin{array}{cc}
\rho\cos\alpha & \sigma\cos\beta \\
\rho\sin\alpha & \sigma\sin\beta
\end{array}\right] =
\left[\begin{array}{cc}
3 & 1 \\
1 & 4
\end{array}\right].
\]
Then we still have that $\sin(\beta-\alpha)= \|h\|/\|v_2\|$ and
\begin{align*}
a &= \|v_1\| \|h\| = \|v_1\| \|v_2\| \sin(\beta-\alpha),\\
 &= \rho\sigma (\cos\alpha \sin\beta - \cos\beta\sin\alpha), \\
 &= \rho\cos\alpha \,\sigma\sin\beta - \sigma\cos\beta\,\rho\sin\alpha,\\
 &= a_{11} a_{22} - a_{21} a_{12} = 3\cdot 4 - 1\cdot 1 = 11.
\end{align*}
For the students to whom the sum formula
\[
\sin(\beta-\alpha) = \cos\alpha\sin\beta-\cos\beta\sin\alpha
\]
is unknown or feels too much pulled out of the hat, I add an exercise/example when we cover complex numbers:
\begin{align*}
e^{(\beta-\alpha)i} &= \cos(\beta-\alpha) + i \sin(\beta-\alpha),\\
 &= e^{\beta i}e^{-\alpha i} = (\cos\beta + i\sin\beta)(\cos\alpha - i\sin\alpha),\\
 &= (\cos\alpha\cos\beta + \sin\alpha\sin\beta) + i(\cos\alpha\sin\beta - \sin\alpha\cos\beta),
\end{align*}
and thus
\begin{align*}
\cos(\beta-\alpha) &= \text{Re}(e^{(\beta-\alpha)i}) = \cos\alpha\cos\beta + \sin\alpha\sin\beta,\\
\sin(\beta-\alpha) &= \text{Im}(e^{(\beta-\alpha)i}) = \cos\alpha\sin\beta - \sin\alpha\cos\beta.
\end{align*}
If you'd like to follow the senior econometrician's advice and avoid complex numbers altogether, there are also beautiful visual proofs of this result, for example in \cite{Nelsen2000}, p. 43. For the students for whom the Euler formula $e^{ix}=\cos x + i\sin x$ is unknown or feels foreign, I add an exercise when we cover Taylor expansions:

\begin{enumerate}
	\item Let $f(x)=\sin(x)$.  Use the Taylor-formula to show that
	\[
	\sin(x) = \sum_{n=0}^{\infty}  \frac{(-1)^n x^{2n+1}}{(2n+1)!}.
	\]
	One finds	
	\[
	f(x) = x -\frac{x^3}{3!} + \frac{x^5}{5!} - \frac{x^7}{7!} + \frac{x^9}{9!} - \ldots
	\]
	\item Find a corresponding expression for $\cos(x)$.
	
	\[
	f(x) = 1 - \frac{x^2}{2} +\frac{x^4}{4!} - \frac{x^6}{6!} + \frac{x^8}{8!} - \ldots
	\]
	
	\item Use the two results to show that 
	\[
	\exp(ix) = \cos(x) + i \sin(x),
 	\]
 	where $i$ is the imaginary unit satisfying $i^2=-1$.

	One finds	
	\begin{align*}
	\exp(ix) &= 1 +ix - \frac{x^2}{2!} - i\frac{x^3}{3!} + \frac{x^4}{4!} + i\frac{x^5}{5!} - \frac{x^6}{6!} - \ldots \\
	 &= \left( 1 - \frac{x^2}{2} +\frac{x^4}{4!} - \frac{x^6}{6!} + \frac{x^8}{8!} - \ldots \right) \\
	 & \qquad + i \left(  x -\frac{x^3}{3!} + \frac{x^5}{5!} - \frac{x^7}{7!} + \frac{x^9}{9!} - \ldots \right).
	\end{align*}
\end{enumerate}

\subsubsection{Cramer's formula}

Equipped with the geometric intuition of the determinant, developed in the last few pages, the students can now appreciate this wonderful visualization, which I found in \cite{Nelsen2000}.

\begin{center}
\begin{tikzpicture}

\filldraw[color=lightgray] (0,0) -- (2,0) -- (3.5,.6) -- (1.5,.6) --
(0,0);

\draw[->,thick] (0,0) -- (2,0) node[anchor=north] {$\scriptstyle a_1$};

\draw[->,thick] (0,0) -- (1.5,.6) node[anchor=north west] {$\scriptstyle
a_2$};

\draw[->,thick] (0,0) -- (.5,1) node[anchor=north west] {$\scriptstyle
a_3$};

\draw[->,thick] (0,0) -- (1.5,3) node[anchor=north west] {$\scriptstyle
x_3 a_3$};

\draw[->] (0,0) -- (0,5);

\draw (2,0) -- (3.5,3);

\draw (1.5,.6) -- (3,3.6);

\draw (2,0) -- (3.5,.6);

\draw (1.5,.6) -- (3.5,.6);

\draw (3.5,.6) -- (5,3.6);

\draw (.5,1) -- (2,1.6) -- (4,1.6) -- (2.5,1) -- (.5,1);

\draw (-1,-1) -- (0,0);

\draw (-1,2) -- (0,3);

\draw (-.8,-.8) -- (-.8,2.2);

\draw (-.8,.7) node[anchor=east] {$\scriptstyle x_3 h$};

\draw (-.5,.5) -- (0,1);

\draw (-.4,-.4) -- (-.4,.6);

\draw (-.4,.1) node[anchor=east] {$\scriptstyle h$};

\draw (-.1,2.9) -- (-.1,2.7) -- (0,2.8);

\draw (-.1,.9) -- (-.1,.7) -- (0,.8);

\draw[dashed] (0,3) -- (1.5,3);

\draw[dashed] (0,1) -- (.5,1);

\draw[dashed] (3,3.6) -- (5.5,4.6) -- (7,4.6);

\draw[dashed] (7.5,4) -- (5,3) -- (3.5,3);

\draw[->,thick] (0,0) -- (5.5,4) node[anchor=south east] {$\scriptstyle
b$};

\draw (1.5,.6) -- (7,4.6);

\draw (2,0) -- (7.5,4);

\draw (3.5,.6) -- (9,4.6);

\filldraw[color=lightgray] (5.5,4) -- (7.5,4) -- (9,4.6) -- (7,4.6) --
(5.5,4);

\draw (5.5,4) -- (7.5,4) -- (9,4.6) -- (7,4.6) -- (5.5,4);

\filldraw[color=lightgray] (1.5,3) -- (3,3.6) -- (5,3.6) -- (3.5,3) --
(1.5,3);

\draw (1.5,3) -- (3,3.6) -- (5,3.6) -- (3.5,3) -- (1.5,3);

\draw (6.15,1) node {$\scriptstyle b = x_1 a_1 + x_2 a_2 + x_3 a_3$};

\draw (8,.5) node {$\scriptstyle \det(a_1,\, a_2,\, b) = \det(a_1,\,
a_2,\, x_3 a_3) = x_3 \det(a_1,\, a_2,\, a_3)$};

\draw (6,-.1) node {$\scriptstyle x_3 = \frac{\scriptstyle \det(a_1,\,
a_2,\, b)}{\scriptstyle \det(a_1,\, a_2,\, a_3)}$};

\end{tikzpicture}\hfill\\
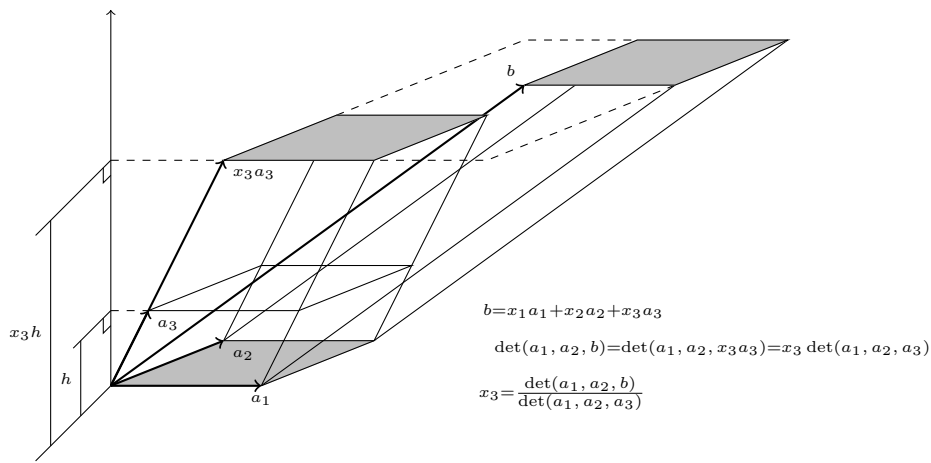
{\tiny From Nelsen (2000), ``Proofs without words.''}
\end{center}

\subsubsection{Polynomials as determinants/volumes}

This example shows that polynomials are at the same time determinants of a set of matrices, and their function values are therefore volumes of parallelepipeds. The students will have to work with the characteristic polynomial in solving the eigenvalue problem, and this example establishes the close link between matrices, determinants, and polynomials. Consider the matrix
{\tiny
\[
A = \left[\begin{array}{cccccc}
a_{n-1}+x & a_{n-2} & a_{n-3} & \cdots & a_1 & a_0 \\
-1 & x & 0 & \cdots & 0 & 0 \\
0 & -1 & x & \cdots & 0 & 0 \\
\vdots & \vdots & \vdots & \ddots & \vdots & \vdots \\
0 & 0 & 0 & \cdots & x & 0 \\
0 & 0 & 0 & \cdots & -1 & x
\end{array}\right] = x I - \left[\begin{array}{cccccc}
-a_{n-1} & -a_{n-2} & -a_{n-3} & \cdots & -a_1 & -a_0 \\
1 & 0 & 0 & \cdots & 0 & 0 \\
0 & 1 & 0 & \cdots & 0 & 0 \\
\vdots & \vdots & \vdots & \ddots & \vdots & \vdots \\
0 & 0 & 0 & \cdots & 0 & 0 \\
0 & 0 & 0 & \cdots & 1 & 0
\end{array}\right].
\]
}
Expand the determinant along the first row
{\tiny
\begin{align*}
\det A &= (a_{n-1} + x) \det\left[\begin{array}{ccccc}
x  & 0 & \cdots & 0 & 0 \\
-1 & x & \cdots & 0 & 0 \\
\vdots & \vdots & \ddots & \vdots & \vdots \\
0  & 0 & \cdots & x & 0 \\
0  & 0 & \cdots & -1 & x \\
\end{array}\right]
 -a_{n-2} \det\left[\begin{array}{ccccc}
-1 & 0 & \cdots & 0 & 0 \\
0  & x & \cdots & 0 & 0 \\
\vdots & \vdots & \ddots & \vdots & \vdots \\
0  & 0 & \cdots & x & 0 \\
0  & 0 & \cdots & -1 & x \\
\end{array}\right]\\
 &+a_{n-3} \det\left[\begin{array}{ccccc}
-1 & x & \cdots & 0 & 0 \\
0 & -1 & \cdots & 0 & 0 \\
\vdots & \vdots & \ddots & \vdots & \vdots \\
0  & 0 & \cdots & x & 0 \\
0  & 0 & \cdots & -1 & x \\
\end{array}\right] + \ldots + (-1)^{n+1} a_0 \det\left[\begin{array}{ccccc}
-1 & x & \cdots & 0 & 0 \\
0 & -1 & \cdots & 0 & 0 \\
\vdots & \vdots & \ddots & \vdots & \vdots \\
0  & 0 & \cdots & -1 & x \\
0  & 0 & \cdots & 0  & -1 \\
\end{array}\right].\\
\end{align*}
}
The right-hand side determinants are easily evaluated since they
reduce to determinants of upper and lower triangular matrices.
\begin{align*}
\det A &= x^n + a_{n-1} x^{n-1} - a_{n-2} (-1) x^{n-2} + \ldots + (-1)^{2n} a_0,\\
 &= x^n + a_{n-1} x^{n-1} + a_{n-2} x^{n-2} + \ldots + a_1 x + a_0.
\end{align*}
Therefore, any $n$-th degree polynomial is the volume of the parallepiped
spanned by the columns of the matrix $A(x)$. The GeoGebra object \url{https://www.geogebra.org/3d/n2ppshze} shows the case $n=3$ with $a_2=1$, $a_1=1$, and $a_0=1$ and the column vectors for the argument values $x=-1$ and $x=3$.

\subsection{Eigenvalues and Lagrange Multipliers}

\subsubsection{The gradient of a quadratic symmetric form}
This example or exercise shows how to apply the definition of differentiability of a function from $\R^n$ to $\R^m$ to symmetric bilinear forms. I find it particularly intuitive because of the close resemblance of the gradient to the univariate situation.  Let $A\in\R^{n\times n}$ be a symmetric matrix. Consider the quadratic symmetric form
\[
q_A(x) = x^T A x = \sum_{i=1}^n\sum_{j=1}^n a_{ij} x_i x_j = \langle x,\, Ax \rangle.
\]
Then,
\begin{align*}
q_A(x+h) &= \langle x+h,\, A(x+h) \rangle = \langle x+h,\, Ax + Ah \rangle\\
 &= \langle x,\, Ax\rangle + \langle x,\, Ah\rangle + \langle h,\, Ax\rangle + \langle h,\, Ah\rangle,\\
\langle x,\, Ah\rangle &= x^T Ah = (x^T Ah)^T = h^T A x = \langle h,\, Ax\rangle .
\end{align*}
Therefore,
\begin{align*}
q_A(x+h) &= \langle x,\, Ax\rangle + 2\langle h,\, Ax\rangle + \langle h\, Ah\rangle,\\
 &= q_A(x) + \langle 2Ax,\, h\rangle + \langle h,\, Ah\rangle ,
\end{align*}
where
\[
\lim_{h \to 0}\frac{\langle h,\, Ah\rangle}{\| h\|} = 0.
\]
The definition of differentiability of a function $f:\R^n\to\R^m$ is that there is a representation
\[
f(x+h) = f(x) + Ch + R(h), \; x,h\in\R^n, C\in\R^{m\times n}, \lim_{h\to 0} R(h)/\|h\|=0.
\] 
The result is that quadratic symmetric forms $q_A: \R^n\to\R$ are totally differentiable with
\[
C^T = \nabla q_A(x) = 2 A x \in \R^{n\times 1}.
\]
Loosely speaking, the quadratic form $x^T A x$ is the multivariate version of the function $f(x)=ax^2$, and so this result is not surprising.

\subsubsection{Extrema of a quadratic symmetric form on the unit sphere}

This example extends the insight into symmetric bilinear forms from the derivation above to an optimization problem. It is an example where Lagrange multipliers turn out to be eigenvalues.  Let $A\in\R^{n\times n}$ be a symmetric matrix, $\lambda\in\R$. Find the extrema of $q_A(x) = x^T Ax$ on the unit sphere
    \[
        S = \{x\in\R^n \mbox{ such that } \|x\|^2 = \langle x,\, x\rangle = 1\}.
    \]
    Since $q_A$ is a continuous function and $S$ a compact set, extrema exist. Write the constraint as
    \[
        g(x) = \langle x,\, x\rangle - 1 = \sum_{i=1}^n x_i^2 - 1 = 0.
    \]
    The partial derivatives and the gradient of the constraint are
    \[
        \frac{\partial g}{\partial x_i} = 2x_i,\qquad \nabla g(x) = 2x.
    \]
   According to the Lagrange multiplier theorem, the necessary condition
    is
    \begin{align*}
    \nabla q_A &= \lambda \, \nabla g,\\
    2Ax &= \lambda 2x,
    \end{align*}
    or
    \[
    \displaystyle Ax = \lambda\, x,
    \]
    which we recognize as the eigenvalue problem. Denote the minimizer of
    $q_A$ on $S$ by $x_0$ and the maximizer by $x_1$.  Then,
    \begin{align*}
    q_A(x_0) &= \langle x_0,\, A x_0\rangle = \langle x_0,\, \lambda_0 x_0 \rangle = \lambda_0\, \langle x_0,\, x_0 \rangle\\
             &= \lambda_0 \mbox{ smallest eigenvalue},\\
    q_A(x_1) &= \langle x_1,\, A x_1\rangle = \langle x_1,\, \lambda_1 x_1 \rangle = \lambda_1\, \langle x_1,\, x_1 \rangle\\
             &= \lambda_1 \mbox{ largest eigenvalue}.
    \end{align*}
    It follows that the smallest and the largest eigenvalue of $A$ are the minimum and the maximum of $q_A(x) = x^T Ax$ on the unit sphere $S$ and at the same time the Lagrange multipliers of this constrained optimization problem.

\subsection{The airborne fraction of anthropogenic CO$_2$ emissions}

In the following, I develop an example of a first-order linear differential equation with a constant homogeneous coefficient and time-dependent inhomogeneous term. A common assumption in climate science is that the terrestrial biosphere (plants) and the ocean absorb atmospheric CO$_2$ linearly,
\begin{align*}
F_{oc}(t) &= \frac{1}{\tau_{oc}} C(t), \\
F_{ld}(t) &= \frac{1}{\tau_{ld}} C(t),
\end{align*}
where $F_{oc,ld}$ denotes absorption by ocean and land ($F$-flux, $oc$-ocean, $ld$-land), and the positive numbers $\tau_{oc}$ and $\tau_{ld}$ denote the average time it takes for land and ocean to absorb one unit of carbon $C$ of atmospheric CO$_2$. 

Under the assumption that anthropogenic CO$_2$ emissions are sufficiently well-described by an exponential function, which is roughly true on the total period from the beginning of industrialization to today, we have,\footnote{On the more recent sample since 1959, a case can be made that emissions are linear, see \cite{Bennedsen2023a}.}
\[
f(t) = f_0 e^{dt}, \; d>0,
\]
it can be shown that the so-called ``\textit{airborne fraction}'' 
\[
AF = \frac{1}{f(t)}\frac{dC}{dt},
\]
or the ratio of anthropogenic CO$_2$ emissions $f(t)$ that is not absorbed by land or ocean and thus remains in the atmosphere (``airborne''), converges to a constant \citep{Bennedsen2019,Bennedsen2023b,Bennedsen2024}.

To approach this claim with the theory of differential equations, we note the global carbon budget equation, which states that anthropogenic CO$_2$ emissions are either absorbed by land, ocean, or the atmosphere. The instantaneous change in atmospheric CO$_2$ is thus given by 
\[
\frac{d C(t)}{dt} = f(t) - \left(\frac{1}{\tau_{oc}} + \frac{1}{\tau_{ld}}\right) C(t).
\]
With $x(t):= C(t)$ and $c:=\tau_{oc}^{-1}+\tau_{ld}^{-1}>0$ obtains
\[
\dot x(t) + c x(t) = f(t).
\]
Multiplying with the integrating factor $e^{ct}$ yields
\[
\dot x(t) e^{ct} + c e^{ct} x(t) = f(t) e^{ct} = \frac{d}{dt} (xe^{ct}).
\]
Integrating, we obtain
\[
xe^{ct} = x_0 + \int_0^t f(s) e^{cs} ds,
\]
or,
\[
x = x_0 e^{-ct} + e^{-ct}\int_0^t f(s) e^{cs} ds.
\]
Plugging in the assumption for emissions that $f(t) = f_0 e^{dt}$, we get
\begin{align*}
x &= x_0 e^{-ct} + e^{-ct}f_0\int_0^t e^{(c+d)s} ds \\
 &= x_0 e^{-ct} + f_0 e^{-ct} \left[ \frac{1}{c+d} e^{(c+d)s} \right]_0^t \\
 &= x_0 e^{-ct} + f_0 e^{-ct} \left( \frac{1}{c+d} e^{(c+d)t} -  \frac{1}{c+d}\right).
\end{align*}
In order to calculate the airborne fraction, we need a solution for $x/f$:
\begin{align*}
\frac{x(t)}{f(t)} &= f_0^{-1} e^{-dt} \left(x_0 e^{-ct} + f_0 e^{-ct} \left( \frac{1}{c+d} e^{(c+d)t} -  \frac{1}{c+d}\right)\right) \\
 &= \frac{x_0}{f_0} e^{-(c+d)t} + \frac{1}{c+d} - \frac{1}{c+d} e^{-(c+d)t} \\
 &= \left(\frac{x_0}{f_0} - \frac{1}{c+d}\right)e^{-(c+d)t} + \frac{1}{c+d}.
\end{align*}
We arrive at
\begin{align*}
AF=\frac{\dot x(t)}{f(t)} &= \frac{f(t) - c x(t)}{f(t)} = 1 - c \frac{x(t)}{f(t)}\\
 &= 1 - c \left(\frac{x_0}{f_0} - \frac{1}{c+d}\right)e^{-(c+d)t} -  \frac{c}{c+d}.
\end{align*}
As $t\to\infty$, the airborne fraction converges to
\[
\lim_{t\to\infty} AF= \lim_{t\to\infty}\frac{\dot x(t)}{f(t)} = 1-\frac{c}{c+d} = \frac{d}{c+d}.
\]

\subsection{The Ramsey model as a recurring example throughout a course}

I use the Ramsey model as a recurring example throughout the third-semester undergraduate course based on \cite{Sydsaeter2008}, since it touches on many elements that we cover. The sequence of examples from the Ramsey model culminates in the optimization problem for the household, for which control theory is necessary. The entire course also culminates, in many respects, with the chapter on control theory, which draws on elements from linear algebra, multivariate calculus, static optimization, and differential equations.

\subsubsection{CRR or CIES utility functions}

This is a preliminary exercise that provides an insight that is needed later to solve the Ramsey model, but of course it is also of independent interest. Solving the exercise requires only methods from first-order ordinary differential equations. It is also the only example of a differential equation where the indicator variable is not time in my course.

Let $U : x \mapsto U(x)\in\R$ be a twice continuously differentiable utility function.  The elasticity of the first derivative $U'$ with respect to the good $x$,
\[
    -\epsilon_x(U') = -\left.\frac{dU'(x)}{U'(x)}\right/\frac{dx}{x} = -\frac{dU'(x)}{dx}\, \frac{x}{U'(x)} = -\frac{U''(x) x}{U'(x)}
\]
    is called the \textit{Arrow-Pratt measure} of relative risk
    aversion.  By requiring constant relative risk aversion (CRR), or constant (intertemporal) elasticity of substitution (C(I)ES), one arrives at a second-order differential equation for $U$:
    \[
        -\frac{U''(x) x}{U'(x)} \equiv \theta \; \Longleftrightarrow\; U''(x) = -\frac{\theta}{x}\, U'(x),\, \theta>0.
    \]
    Substitute $v(x) := U'(x)$, then,
    \[
        v'(x) = -\frac{\theta}{x}\, v(x)
    \]
    Solve this first-order DE for $v$ and then integrate $v$ to obtain $U$.  Distinguish the cases $\theta\neq 1$ and $\theta =1$.  The resulting class of functions is the family of CRR or C(I)ES utility functions.  The higher $\theta$, the faster the relative decline in $U'(x)$ in response to increases in $x$.  That is, increases in $x$ are less welcome than with a lower value of $\theta$.
    
We obtain
    \begin{align*}
        v'(x) &= -\frac{\theta}{x}\, v(x),\\
        \int \frac{v'(x)}{v(x)}\, dx &= \int \frac{1}{v(x)}\, dv = -\int \frac{\theta}{x} dx,\\
        \log v(x) &= - \theta(\log x + c),\, c\in\R,\\
        v(x) &= e^{-\theta c}\, x^{-\theta} =: k_0 x^{-\theta}.
    \end{align*}
    Then,
    \begin{align*}
        U(x) &= \int U'(x) \, dx = \int v(x)\, dx\\
         &= \int k_0 x^{-\theta}\, dx = \frac{x^{1-\theta}}{1-\theta}k_0  + k_1,\, k_1\in\R,
    \end{align*}
    for $\theta\neq 1$.  For $\theta =1$,
    \[
        \int k_0 x^{-1}\, dx = k_0 \log x + k_1.
    \]
    This is the family of CRR or C(I)ES utility functions.  The higher $\theta$, the faster the relative decline in $U'(x)$ in response to increases in $x$.  That is, increases in $x$ are less welcome than with a lower value of $\theta$.

\subsubsection{Linearization around steady-state as example for Taylor expansion}

The students' first point of contact with the Ramsey model is in multivariate calculus when we cover the Taylor expansion. I introduce the system of differential equations that describes capital and consumption dynamics in the model, even though the students have not seen differential equations at this point. Since the task is to understand the linearization of the right-hand sides of the differential equations around the steady state and not to solve the equations, this is not an issue. I usually spend a few words explaining the idea behind differential equations and that we are going to come back to this later in the course. I recurrently use Cobb-Douglas production and CES utility functions as examples throughout the course, so the students are familiar with them at this point. They also take Macroeconomics I in the same semester. I use three slides to introduce the following material.

In macroeconomics, the Ramsey model is a key model to describe consumption and production in a stylized economy \citep{Ramsey1928,Koopmans1963,Cass1965}, with textbook treatments, for example, in \cite{Barro2004} and \cite{Romer2006}.

In equilibrium, the economy is described by a system of non-linear ordinary differential equations for capital per unit of effective labor $k(t)$ and consumption per unit of effective labor $c(t)$, as functions of time.
\begin{align*}
\frac{d}{dt} \log k(t) &= A k^{\alpha-1} - \frac{c}{k} - (\delta + \alpha_L + \alpha_T),\\
\frac{d}{dt} \log c(t) &= \frac{1}{\theta} \left( \alpha A k^{\alpha-1} - (\delta + \rho + \theta \alpha_T) \right).
\end{align*}
\begin{itemize}
\item $A,\, \alpha>0$ are the parameters of a Cobb-Douglas production function, 
\item $\theta>0$ is the parameter of a CES utility function,
\item $\delta >0$ is the depreciation of the capital stock,
\item $\alpha_L>0$ is the growth rate of the labor force,
\item $\alpha_T>0$ is the rate of technological progress, 
\item $\rho>0$ is the rate of time preference.
\end{itemize}

The \textit{steady state} $(k(t),\, c(t)) \equiv (k^\ast,\, c^\ast)$, where neither capital nor consumption change through time, is defined by
\[
\frac{d}{dt} \log k(t) = \frac{d}{dt} \log c(t) = 0.
\]
The fact that the system of differential equations is non-linear is often cumbersome, for example in numerical analyses.  For many purposes, the system is linearized around the steady state.  That is, the function
\begin{align*}
f : \R^2 &\longrightarrow \R^2,\\
\left[ \begin{array}{c} \log k \\ \log c \end{array} \right] &\longmapsto \left[ \begin{array}{c} \frac{d}{dt} \log k \\ \frac{d}{dt} \log c \end{array} \right],
\end{align*}
is expanded around $(\log k^\ast,\, \log c^\ast)$ to first order. 

The resulting first-order Taylor expansion is a linear system of differential equations:
\[
\left[ \begin{array}{c} \frac{d}{dt} \log k \\ \frac{d}{dt} \log c \end{array} \right] = \left[ \begin{array}{cc} \rho-\alpha_L-(1-\theta)\alpha_T & \delta+\alpha_L+\alpha_T-\frac{1}{\alpha}(\delta+\rho+\theta\alpha_T) \\ \frac{\alpha-1}{\theta}(\delta+\rho+\theta\alpha_T) & 0 \end{array} \right] \left[ \begin{array}{c} \log \frac{k}{k^\ast} \\ \log \frac{c}{c^\ast} \end{array} \right].
\]
To fill in the blanks in this example is the following exercise in a homework set. 

\begin{enumerate}
	\item Show that the Cobb-Douglas production function $Y=AK^\alpha L^{1-\alpha}$ can be written in units of effective labor
	\[
	y = A k^\alpha,
	\]
	where $y=Y/L$ and $k=K/L$.
	\item Show that the steady-state values for capital and consumption in the Ramsey model are given by 
	\begin{align*}
	k^\ast &= \left[ \frac{1}{\alpha A} (\delta + \rho + \theta\alpha_T) \right]^{\frac{1}{\alpha-1}},\\
	c^\ast &= A (k^\ast)^\alpha - (\delta + \alpha_L + \alpha_T) k^\ast.
	\end{align*}
	
	Set $\frac{d}{dt} \log k(t) = \frac{d}{dt}\log c(t) = 0$ and solve for $c$ and $k$.
	
	\item Show that the Jacobian of
	\begin{align*}
	f : \R^2 &\longrightarrow \R^2,\\
	\left[ \begin{array}{c} \log k \\ \log c \end{array} \right] &\longmapsto \left[ \begin{array}{c} Ak^{\alpha-1} - \frac{c}{k} - (\delta + \alpha_L + \alpha_T) \\ \frac{1}{\theta} \left(\alpha A k^{\alpha-1} - (\delta + \rho + \theta\alpha_T)\right) \end{array} \right],
	\end{align*}
	evaluated in $(\log k^\ast,\, \log c^\ast)$ is given by
	\[
	J_f(\log k^*,\,\log c^* ) = \left[ \begin{array}{cc} \rho - \alpha_L - (1-\theta) \alpha_T & \delta + \alpha_L + \alpha_T - \frac{1}{\alpha} (\delta + \rho + \theta\alpha_T) \\ \frac{\alpha-1}{\theta} (\delta + \rho + \theta\alpha_T)  &   0 \end{array} \right].
	\]
	\item Show that
	\[
	f(\log k^*,\, \log c^*) = (0,\, 0)^T.
	\]
\end{enumerate}

Hint: 	Rewrite the function $f$ as
\[
\left[\begin{array}{c}
\log k \\ \log c
\end{array}\right]
\mapsto
\left[\begin{array}{c}
A e^{(\alpha-1)\log k} - e^{\log \frac{c}{k}} - (\delta+\alpha_L+\alpha_T) \\ \frac{1}{\theta}\left(\alpha A e^{(\alpha-1)\log k} - (\delta + \rho + \theta\alpha_T )\right)
\end{array}\right]
\]
with Jacobian
\[
J_f = \left[\begin{array}{cc}
A(\alpha-1)k^{\alpha-1} + \frac{c}{k} & -\frac{c}{k} \\[.4ex]
\frac{\alpha(\alpha-1)}{\theta} A k^{\alpha -1} & 0
\end{array}\right].
\]

\subsubsection{Understanding the household's optimization problem}

When we cover first-order linear ordinary differential equations (Chapter 5 in \cite{Sydsaeter2008}), I spend three slides on explaining the household's optimization problem. The students cannot solve it yet, but they can understand the objects that figure in the problem.

In the Ramsey model, define the following objects.
\smallskip

\begin{tabular}{ll}
$\mathsf{L}(t) = e^{\alpha_L t}$ & work force per representative household\\
$\alpha_L>0$ & population growth rate \\
$C(t)$ & aggregate consumption of the representative household\\
$\texttt{c}(t)=\frac{C(t)}{\mathsf{L}(t)}$ & per-capita consumption\\
$u(x)$ & CRRA, CES utility function  \\
$\rho>0$ & rate of time preference \\
$A(t)$ & total assets of the representative household\\
$\text{a}(t) = \frac{A(t)}{\mathsf{L}(t)}$ & per-capita assets \\
$r(t)$ & return rate on assets \\
$\texttt{w}(t)$ & wage rate per unit of effective labor
\end{tabular}

Assets can be ownership of capital or loans to firms; both earn the return rate $r(t)$.  Negative values of $\texttt{a}$ are possible; they represent debt to firms. The instantaneous utility function is given by
\[
u(x) = \frac{x^{1-\theta}-1}{1-\theta}.
\]
The utility of the representative household is then given by integrating over the individual utility multiplied by the number of individuals and a discount factor to account for time preference.  Assuming the planning horizon $[t_0,\, t]$, the objective function of the representative household is
\[
U(t_0,\, t,\, a,\, \texttt{c}) = \int_{t_0}^t u(\texttt{c}(s)) \mathsf{L}(s) e^{-\rho s} ds.
\]
This objective function is maximized subject to a constraint:  The present value of all consumption during the planning interval plus the present value of any assets left at time $t$ must equal the present value of all income, plus the assets owned at the beginning of the period, $t_0$.

The instantaneous change in the total asset position $A(t)$ of the representative household is given by returns on the total asset position, $r(t)A(t)$, plus total wages per household, $\texttt{w}(t)\mathsf{L}(t)$, minus household consumption $C(t)$:
\[
\dot{A}(t) = r(t)A(t) + \texttt{w}(t)\mathsf{L}(t) - C(t),
\]
or in per-capita terms,
\[
\dot{\texttt{a}}(t) = r(t)\texttt{a}(t) + \texttt{w}(t) - \texttt{c}(t) - \alpha_L \texttt{a}(t).
\]
(Show that $\dot{\texttt{a}}(t) = \dot{A}(t) - \alpha_L \texttt{a}(t)$.) The solution for $\texttt{a}(t)$ is obtained by variation of constants as
\[
\texttt{a}(t) = \texttt{a}(t_0)e^{\int_{t_0}^t (r(s)-\alpha_L)ds} + e^{\int_{t_0}^t (r(s)-\alpha_L)ds} \int_{t_0}^t (\texttt{w}(s)-\texttt{c}(s))\, e^{-\int_{t_0}^s (r(\tau)-\alpha_L)d\tau}\, ds,
\]
or
\[
\texttt{a}(t)e^{-\int_{t_0}^t (r(s)-\alpha_L)ds} + \int_{t_0}^t \texttt{c}(s) e^{-\int_{t_0}^s (r(\tau)-\alpha_L)d\tau}\, ds = \texttt{a}(t_0) + \int_{t_0}^t \texttt{w}(s) e^{-\int_{t_0}^s (r(\tau)-\alpha_L)d\tau}\, ds.
\]
The last equation expresses the constraint described verbally above.

\subsubsection{Solving the linearized system of differential equations}

When we get to second-order linear differential equations (Chapter 6 in \cite{Sydsaeter2008}), we have the tools to solve the linearized system we found in the first exercise:

In the Ramsey model, we found a linearization of the logarithmic capital stock $\log k(t)$ and the logarithmic consumption $\log c(t)$ around their steady-state values $(\log k^*$, $\log c^*)$ as
\[
\dot{x}(t) = A x(t) + b,
\]
where
\[
x(t) = \left[\begin{array}{c}\log k(t) \\ \log c(t)\end{array}\right],
\]
and
\[
A = \left[\begin{array}{cc}
a_{11} & a_{12} \\
a_{21} & a_{22}
\end{array}\right],
\]
with entries
\begin{align*}
a_{11} &= \rho - \alpha_L - (1-\theta) \alpha_T,\\
a_{12} &= \delta + \alpha_L + \alpha_T - \frac{1}{\alpha}(\delta + \rho + \theta\alpha_T),\\
a_{21} &= \frac{\alpha -1}{\theta} (\delta + \rho + \theta\alpha_T),\\
a_{22} &= 0.
\end{align*}
The inhomogeneous term is
\[
b = -A \left[\begin{array}{c} \log k^* \\ \log c^* \end{array}\right],
\]
where
\begin{align*}
k^* &= \left[  \frac{1}{\alpha A} (\delta + \rho + \theta\alpha_T)  \right]^{\frac{1}{\alpha -1}},\\
c^* &= A (k^*)^\alpha - (\delta + \alpha_L + \alpha_T) k^*.
\end{align*}

\begin{quote}
\textbf{Exercise:} \textit{Determine the eigenvalues $\lambda_{1}$, $\lambda_2$ and eigenvectors $v_{1}$, $v_2$ of $A$ (as functions of the $a_{ij}$).  Show that the matrix $A$ is diagonalizable if and only if}
\[
a_{11}^2 \neq -4 a_{12}a_{21}.
\]
\end{quote}

The solutions of the linearized system have the form
\[
[v_1 e^{\lambda_1 t},\, v_2 e^{\lambda_2 t}  ].
\]
The inhomogeneous equation 
\[
\dot{x}(t) = A x(t) + b,
\]
where $b = - A x^*$ and $x^*=(\log k^*,\, \log c^*)$, is solved with the substitution $z= \log k -\log k^*$ and $w= \log c -\log c^*$.

Eigenvalues and -vectors of $A$ are given by
\begin{align*}
\lambda_1 &= \frac{a_{11}}{2} + \sqrt{\frac{a_{11}^2}{4} + a_{12}a_{21}}, & v_1 &= \left[\begin{array}{c} -a_{12} \\ \frac{a_{11}}{2} - \sqrt{\frac{a_{11}^2}{4} + a_{12}a_{21}} \end{array}\right],\\
\lambda_2 &= \frac{a_{11}}{2} - \sqrt{\frac{a_{11}^2}{4} + a_{12}a_{21}}, & v_2 &= \left[\begin{array}{c} -a_{12} \\ \frac{a_{11}}{2} + \sqrt{\frac{a_{11}^2}{4} + a_{12}a_{21}} \end{array}\right].
\end{align*}

\subsubsection{Solving the optimization problems and closing the model}

When we have covered dynamic optimization problems with an infinite horizon in the chapter on control theory (Chapter 9 in \cite{Sydsaeter2008}), I present the students with this ``soft'' exercise. In the oral exam, the students draw a card with a subject out of a dozen or so, covering different parts of the course. They prepare an 8-minute presentation on each subject, so that they have a degree of control.  I agree in advance with the students that I will only talk about the Ramsey model if they choose it in their presentation when they draw a subject in the chapter on systems of differential equations or in control theory. In the ten years that I have been teaching this course, no student has chosen to do so, even though some have expressed appreciation of the coverage of the model.

\smallskip

This exercise involves understanding the text and following the steps to arrive at the displayed equations. Throughout this class, we have used elements of the Ramsey model in macroeconomics as motivation.  We are now in a position to derive the first-order conditions for optimality in this model, which are a system of differential equations. The notation used here resembles the one in Barro and Sala-i-Martin (2004), Economic Growth, 2nd ed.    

\smallskip

The general structure of the model is that firms and households (more precisely, the representative firm and the representative household) have independent optimization problems. Firms maximize profits with respect to capital input.  Profits equal price-times-output minus costs, which are returns on capital to be paid to the owner of the capital, capital depreciation, and wages for labor employed.  Output is a (Cobb-Douglas production) function of capital and labor.  Firms take the price for the final product as given, and it is customarily normalized to one.  They also take the return rate on capital and the wage rate per unit of labor as given, meaning that they cannot choose these values.  It turns out that the firm's optimization problem can be solved by means of basic calculus and need not be cast as an optimal control problem.

\smallskip

Households maximize utility for their working members over an infinite time horizon. Utility is a function of consumption, which can be paid for by either wages for labor or by returns on assets that the household owns.  Revenue from either source that is not consumed is saved, builds up assets, and thus becomes the source of future returns.  We have studied the dynamics of the asset build-up as an example of a univariate inhomogeneous linear ordinary differential equation and determined the solution in an example in the chapter on first-order ordinary differential equations.  The household's optimization problem is an optimal control problem.  The objective function depends on the state variable assets, which in turn depends on the control variable consumption.  The differential transition equation for the budget constraint in the example in the chapter on first-order differential equations describes the relation between assets and consumption and constitutes the main constraint of the optimization problem.  The household takes wages and returns on assets as given.

\smallskip

The two independent optimization problems of the firms and of the households therefore take two variables as given:  returns on assets and wages in the case of the household, and returns on capital and wages in the case of the firm.  These loose ends are now tied down by imposing equilibrium.  The equilibrium conditions ``close'' the model. The wage variable is obviously the same for both entities, firms and households, and by principles of microeconomics, in equilibrium it must equal the marginal product of labor, which can be determined from the production function.  Returns on assets and returns on capital are imposed to coincide by requiring that in equilibrium, the total asset stock equals the total capital stock. This models ownership of the productive capital stock by the households.  Note that the labor force is not a choice variable in the Ramsey model.  Economists say that labor force is determined exogenously, meaning that it is not chosen by any agent inside the model. The labor force simply grows exponentially at a constant rate.  Exogenous parameters that describe agent's behavior, but are not consciously and optimally chosen by the agent, such as the rate of time preference or the parameter of risk aversion in the utility function, are often called ``deep parameters.''

\smallskip

In the example on the household budget constraint in the Ramsey model in the chapter on first-order differential equations, we introduced the following objects:

\smallskip

\begin{tabular}{ll}
$\mathsf{L}(t) = e^{\alpha_L t}$ & work force per representative household\\
$\alpha_L>0$ & population growth rate \\
$C(t)$ & aggregate consumption of the representative household\\
$\mathtt{c}(t)=\frac{C(t)}{\mathsf{L}(t)}$ & per-capita consumption\\
$u(x)$ & $=(x^{1-\theta}-1)/(1-\theta)$, CES utility function \\
$\rho>0$ & rate of time preference \\
$A(t)$ & total assets of the representative household\\
$\mathtt{a}(t) = \frac{A(t)}{\mathsf{L}(t)}$ & per-capita assets \\
$r(t)$ & returns on assets \\
$\mathtt{w}(t)$ & wage rate per units of labor
\end{tabular}

\smallskip

We will need the following additional objects; most of them are related to the firm's optimization problem.

\smallskip

\begin{tabular}{ll}
$K(t)$ & capital stock \\
$\delta > 0$ & capital depreciation rate \\
$Y(t)$ & output \\
$\mathtt{T}(t) = e^{\alpha_T t}$ & technology\\
$\alpha_T>0$ & rate of technological progress \\
$L(t)$ & $(=\mathsf{L}(t)\mathtt{T}(t)=e^{(\alpha_L + \alpha_T) t})$ effective labor
\end{tabular}

\smallskip

It is common to work in \textit{units of effective labor}.  We define
\[
y := \frac{Y}{L},\qquad k:= \frac{K}{L},\qquad c:=\frac{C}{L},
\]
as output, capital, and consumption, respectively, in units of effective labor. Note that $L(t) = e^{(\alpha_L + \alpha_T)t} $, and that the initial labor force is normalized to one unit.  Then,
\[
\mathtt{c} = \frac{C}{\mathsf{L}} = \frac{C}{e^{\alpha_L t}},\qquad c = \frac{C}{e^{\alpha_L + \alpha_T} t} = \mathtt{c}e^{-\alpha_T t}.
\]
Analogously, we define $\mathtt{k}=K/\mathsf{L}$, so that $k=\mathtt{k}e^{-\alpha_T t}$.

\smallskip

We begin with the firm's problem.  The production function is of the Cobb-Douglas type, where technology augments labor
\[
Y(t) = F(K(t),\, L(t)) = A L(t)^{1-\alpha} K(t)^{\alpha} = A (\mathsf{L}(t)\mathtt{T}(t))^{1-\alpha} K(t)^{\alpha},
\]
so that the total technological impact on labor productivity is $A(\mathtt{T}(t))^{1-\alpha}$.  In units of effective labor,
\[
y = \frac{Y}{L} = F(k,\, 1) =: f(k) = A k^\alpha.
\]
Firms pay return on capital (to their owners, who are the households), lose depreciation on capital stock, and pay wages.  The price of the final good is normalized to one. We then have the profit function
\begin{align*}
\Pi &= F(K,\, L) - (r+\delta) K - \mathtt{w}\mathsf{L}\\
 &= L f(k) - (r+\delta) Lk - \mathtt{w}Le^{-\alpha_T t}\\
 &= L \left( f(k) - (r+\delta) k - \mathtt{w}\, e^{-\alpha_T t} \right).
\end{align*}
Maximizing this profit function over any time interval is equivalent to maximizing at every single point in time since there is no state variable involved that could transfer residual amounts to the future or borrow from the future.  Differentiation with respect to $k$ and setting the derivative equal to zero results in the first-order condition
\begin{equation}\label{Eq:RamseyFirms1}
r(t) = f'(k(t)) - \delta,
\end{equation}
which determines the return on productive capital as the marginal product minus depreciation.

\smallskip

Turning to the household's problem, the state variable---assets---can shift consumption, and thus utility, to the future or borrow consumption, and thus utility, from the future. (Negative values of $a(t)$ are permitted and interpreted as household debt to firms.)  The household's problem is a proper optimal control problem. The problem is
\begin{align*}
\max_{\mathtt{c}} U(t_0,\, a,\, \mathtt{c}) &= \int_{t_0}^\infty u(\mathtt{c}(s)) \mathsf{L}(s) e^{-\rho s} ds,\\
\mbox{\rm such that} \quad &\dot{\mathtt{a}}(t) = r(t)\mathtt{a}(t) + \mathtt{w}(t) - \mathtt{c}(t) -\alpha_L \mathtt{a}(t),\\
 & \mathtt{a}(0) = \mathtt{a}_0.
\end{align*}
Assets per capita increase due to returns $r$ on assets, the difference between wages and consumption, and the growth rate of the labor force (which enters as the denominator of $\mathtt{a}$).  There is an additional constraint that rules out Ponzi schemes (meaning constant borrowing at a faster rate than the interest rate, so that at the end of the period, a large outstanding debt remains):
\[
\lim_{t\to\infty} \mathtt{a}(t) e^{-\int_{t_0}^t (r(s)-\alpha_L)\, ds} \ge 0.
\]
It will turn out later that this condition is met by the transversality condition.  The Hamiltonian is given by
\[
H(t,\, \mathtt{a},\, \mathtt{c},\, \nu) = u(\mathtt{c}) e^{-(\rho-\alpha_L)t} + \nu \left[ \mathtt{w} + (r-\alpha_L)\mathtt{a} - \mathtt{c} \right].
\]
The first-order conditions according to the maximum principle are:
\begin{align}
\frac{\partial H}{\partial \mathtt{c}} &= e^{-(\rho-\alpha_L)t} u'(\mathtt{c}) - \nu = 0,\label{Eq:RamseyHH1}\\
\dot{\nu} &= -\frac{\partial H}{\partial\mathtt{a}} = -\nu(r-\alpha_L),\label{Eq:RamseyHH2}\\
\lim_{t\to\infty}  \nu(t)\mathtt{a}(t) &= 0.\label{Eq:RamseyHH3}
\end{align}
Differentiate \eqref{Eq:RamseyHH1} with respect to time, plug into \eqref{Eq:RamseyHH2}, and remember from the exercise on constant relative risk aversion that
\[
-\frac{u''(\mathtt{c}) \mathtt{c}}{u'(\mathtt{c})} = \theta,
\]
then
\begin{equation}\label{Eq:RamseyHH4}
\frac{\dot{\mathtt{c}}}{\mathtt{c}} = \frac{1}{\theta}(r(t)-\rho).
\end{equation}
From \eqref{Eq:RamseyHH2}, we have
\[
\nu(t) = \nu(0) e^{-\int_{t_0}^\infty (r(s)-\alpha_L)\, ds},
\]
so that the transversality condition reads
\[
\lim_{t\to\infty} \mathtt{a}(t) e^{-\int_{t_0}^t (r(s)-\alpha_L)\, ds} = 0,
\]
which ensures that the no-Ponzi condition, namely that this expression must be positive or zero. is satisfied. To communicate intuition, I refer informally to the limit as ``the end'' of the planning horizon. Households do not want to have positive assets left over at the end, which is why strictly positive values are ruled out. If positive assets were left at the end in the optimal solution, their shadow price $\nu$ would have to be zero. As long as the shadow price is positive, households do not want to waste assets at the end. In this sense, the transversality condition is the continuous-time and infinite-horizon equivalent of the complementary slack condition in the Kuhn-Tucker conditions.

\smallskip

We now ``close'' the model by imposing equilibrium conditions.  One is the microeconomic insight that in equilibrium, the wage rate $\mathtt{w}(t)$ must equal the marginal product of labor, which we obtain from the production function:
\begin{align*}
Y &= F(K,\, L) = L f(k),\\
\frac{\partial F}{\partial \mathsf{L}} &= \frac{\partial L}{\partial \mathsf{L}} f(k) + L f'(k) \frac{\partial}{\partial\mathsf{L}} \frac{K}{L} = e^{\alpha_T t}(f(k)-f'(k)k) = \mathtt{w}(t).
\end{align*}
The second equilibrium condition states that assets equal capital $\mathtt{a}(t) = \mathtt{k}(t)$ for all $t$, and thus
\[
\dot{\mathtt{k}} = \dot{\mathtt{a}} = (r-\alpha_L) \mathtt{a} + \mathtt{w} - \mathtt{c}.
\]
By definition of $\mathtt{k}$ and $k$, we have
\begin{align}
\dot{k} &= \dot{\mathtt{k}} e^{-\alpha_T t} - \mathtt{k} \alpha_T e^{-\alpha_T t},\notag\\
 &= ((r-\alpha_L) \mathtt{a} + \mathtt{w} - \mathtt{c})e^{-\alpha_T t} - \mathtt{k} \alpha_T e^{-\alpha_T t},\notag\\
 &= f(k) - c - (\delta + \alpha_L + \alpha_T) k,\label{Eq:RamseyK}
\end{align}
where the last equation uses the equilibrium condition for $\mathtt{w}$, the first-order condition \eqref{Eq:RamseyFirms1} from the firm's problem, and $\mathtt{a}=\mathtt{k}$.

\smallskip

From the household's problem, we distilled the condition \eqref{Eq:RamseyHH4}. Verifying that
\[
\frac{\dot{c}}{c} = \frac{\dot{\mathtt{c}}}{\mathtt{c}} - \alpha_T,
\]
and using \eqref{Eq:RamseyFirms1} from the firm's problem, we obtain
\begin{equation}\label{Eq:RamseyC}
\frac{\dot{c}}{c} = \frac{1}{\theta} (f'(k) - \delta - \rho - \theta\alpha_T).
\end{equation}
Inserting the Cobb-Douglas specification $f(k) = ak^\alpha$, Equations \eqref{Eq:RamseyK} and \eqref{Eq:RamseyC} result in the system of non-linear ordinary differential equations stated in the example in our first encounter with the model:
\begin{align*}
\frac{d}{dt} \log k(t) &= \frac{\dot{k}}{k} = A k^{\alpha-1} - \frac{c}{k} - (\delta + \alpha_L + \alpha_T),\\
\frac{d}{dt} \log c(t) &= \frac{\dot{c}}{c} =\frac{1}{\theta} \left( \alpha A k^{\alpha-1} - (\delta + \rho + \theta \alpha_T) \right).
\end{align*}

\section{Conclusion}

In this chapter, I have reflected on the teaching of mathematics for students of economics. I outlined an ideal course sequence if the objective is to enable the student to understand the economic literature at the undergraduate level and contribute to it at the graduate level. I also briefly discussed mathematical modeling in economics. Three years of experience with teaching Mathematics for Economists online was summarized. The second half of the chapter offered a repository of examples and exercises for instructors of Mathematics for Economists, freely available for use. They covered a broad range of subjects in linear algebra, analysis, differential equations, and dynamic optimization. The Ramsey model figured prominently among these examples.

\section*{Acknowledgements}
I would like to thank my students and teaching assistants in 7610 Mathematics for Economists (LSU) and 2622 Mathematics for Economists (AU) over the years. The former TAs Rebekka Jørgensen, Marie Møller Schmidt, and Hannah Segato have read and commented on a draft of this paper. I am grateful for careful reading and comments from Ludwig Arnold, Mikkel Bennedsen, Jan Magnus, and Martin Wagner. The usual disclaimer applies.

\clearpage 

{\small 
\bibliographystyle{chicago}
\bibliography{math_econ_ped}
}

\end{document}